\documentstyle[12pt]{amsart}

\numberwithin{equation}{section}
\theoremstyle{plain}
\newtheorem{Def}{Definition}[section]
\newtheorem{Prop}[Def]{Proposition}
\newtheorem{Lem}[Def]{Lemma}

\newtheorem{Thm}[Def]{Theorem}
\newtheorem{Cor}[Def]{Corollary}

\def\fair{\hat{f}_{\text{Ai}}}
\def\fbess{\hat{f}_{\text{Be}}}
\def\bair{B_{\text{Ai}}}
\def\bbes{B_{\text{Be}}}
\def\too{\longrightarrow}
\def\d{\partial}
\def\C{{\Bbb C}}

\def\grad{Gr^{ad}}

\def\lair{L_{\text{Ai}}}
\def\lbess{L_{\text{Be}}}
\def\vecair{\vec{f}_{\text{Ai}}}
\def\vecbess{\vec{f}_{\text{Be}}}
\def\tr{\text{tr}}
\def\hp{\hat{P}}
\def\hq{\hat{Q}}
\def\tp{\tilde{P}}
\def\tq{\tilde{Q}}
\def\dtr{d\kern 3pt\tr}
\def\kair{K_{0,\text{Ai}}}
\def\psibess{\Psi_{\text{Be},W}}
\def\psiair{\Psi_{\text{Ai},W}}
\def\nabair{\nabla_{\text{Ai}}}
\def\nabbess{\nabla_{\text{Be}}}

\begin{document}

\title[Bispectral Involutions, \today]
{Calogero-Moser pairs and the Airy and Bessel bispectral involutions}
\author{Mitchell Rothstein}
\address{Department of Mathematics\\ University of Georgia\\ Athens, GA 
30602. Temporary Address: Department of Mathematics\\ Harvard University\\
Cambridge, MA\\ 02138}
\email{rothstei@@math.uga.edu}
\thanks{Research supported by NSF Grant No. 58-1353149} 
\date{\today}

\maketitle

\begin{abstract}Explicit formulae are given for the Airy and
Bessel bispectral involutions, in terms of Calogero-Moser pairs.
Hamiltonian structure of the motion of the poles
of the operators
is discussed.
\end{abstract}

\section{Introduction}

This paper  follows upon the study of the Airy bispectral
involution made in \cite{KR}.   There we
gave an analogue, for  arbitrary rank, of the
rank-one bispectral involution developed by Wilson
 \cite{W1}.  Recently, \cite{W2}, Wilson has established a
relationship between the rank-one bispectral involution and the
complex analogue of the Calogero-Moser phase space.  
This relationship leads to  explicit formulae for the Baker function
and the corresponding  involution, 
which make many  important features manifest.  
As shown below, similar results hold for 
bispectral algebras
obtained from generalized Airy and Bessel operators.

Given a positive integer $n$, define
$C_n$ to be the
quotient, under conjugation by $Gl(n,\C)$, of the space of pairs ({\it
Calogero-Moser pairs}) of $n\times n$
complex  matrices, $(P,Q)$, such that
\begin{equation}\label{cm}
\text{rank}([P,Q]-I)=1\ .    
\end{equation}
This is the complex analogue of the definition in \cite{KKS}, in
which $P$ and $Q$ are taken to be hermitian, and $I$ is replaced
by $iI$.  
  Define $\grad$ (see \cite{W1,W2}) to be
the subspace of the Sato grassmannian  (\cite{DJKM,Sa,SW})
corresponding
to Krichever data whose spectral curve is rational and unicursal (no nodes).
The Baker function for such data is always of the form
$e^{xz}p(x,z)$, where $p(x,z)$ is rational and separable. By 
separable we mean its denominator is a product $q(z)\tau(x)$.
Wilson has proved

\begin{Thm}[\cite{W2}, Thm.  5.1]
There is a one-to-one correspondence 
\begin{equation}\label{correspondence}
\cup_{n=0}^{\infty}C_n\leftrightarrow \grad\ ,
\end{equation}
such  that a point $W\in \grad$ corresponds to a point $(P,Q)$ if and only
if its Baker function, $\psi_W$ is given by
\begin{equation}\label{baker formula}
\psi_W(x,z)=e^{xz}\text{Det}(I-(zI-Q)^{-1}(xI-P)^{-1})\ .
\end{equation}
\end{Thm}

Remark:  The ``generic'' rational unicursal curve has simple cusps.
These correspond to $Q$ with distinct eigenvalues.  In the hermitian case this
is automatic. An important aspect of
Wilson's theorem is the observation that in the complex case,
nonsemisimple $Q$ correspond to nonsimple cusps, or in  physical terms,
collision of the Calogero-Moser particles, which are now moving
in $\Bbb C$ rather than $\Bbb R$.

 The  spectral algebra $\cal R_W$ is the algebra 
of differential operators $L(x,\d)$ such that $L\psi_W=f(z)\psi_W$ for
some function $f(z)$. 
Wilson's result 
proves that the spectral algebra of any point $W\in\grad$ is bispectral 
in the sense of \cite{DG}.   That is, the spectral algebra of $\psi_W(z,x)$ is
also nontrivial. 

 It proves much more, for it says
that $\psi_W(z,x)$ is also a Baker function, namely the one corresponding
to $(Q^T,P^T)$.   The involution $(P,Q)\mapsto (Q^T,P^T)$ is clearly
antisymplectic with respect to the symplectic form 
\begin{equation}
\omega=tr(dP\wedge dQ)\ .
\end{equation}
This symplectic structure is an important example of ``unreduction''
\cite{KKS}.    Namely, the
Calogero-Moser hierarchy \cite{AMM} is a completely integrable hamiltonian
system defined on the quotient space $C_n$.   The hamiltonians
are rather complicated in the reduced variables, but on the level
of matrices $(P,Q)$ they are given simply by
the hamiltonians are $h_n=\tr(P^n)$.
 Moreover, the
involution $((P,Q))\too((Q^T,P^T))$ is the linearizing map for the Calogero-Moser
particle system \cite{AMM,KKS}.  
Thus, Wilson's result gives the best proof that
this linearizing map and the bispectral involution are one and the same (cf. 
\cite{Ka}).  As we shall see,  the corresponding involutions,
 in terms of an auxilliary monic polynomial
$\rho(t)$,  are given in the Airy
and Bessel cases respectively by
\begin{align}
(P,Q)&\mapsto (P^T,\rho(P^T)-Q^T)\ ,\\
(P,Q)&\mapsto (
(QP\rho(QP)^{-1} Q)^T,(Q^{-1}\rho(QP))^T)\ .\end{align}

\vskip 12pt
\noindent Acknowledgments:  The author wishes to acknowledge valuable
discussions with J. Harnad, A. Kasman, R. Varley and G. Wilson.

\section{The Airy case}\label{airy case}

The rank of a commutative  algebra
$\cal R$ of ordinary differential operators is 
defined to be the greatest common divisor of the orders of its elements.
The true rank of $\cal R$ is the rank of its centralizer \cite{LP,PW}.   
For instance, fix a monic $r$th order polynomial $\rho(t)$, and define
\begin{equation}
\lair=\rho(\d)-x\ \ \ \ \ \text{(Airy)}\label{airyop}\end{equation}
Then $\Bbb C[\lair]$ has true rank $r$ \cite{KR}.   Moreover, this algebra
is bispectral.  Indeed,   for any $f\in ker(\lair)$, define
\begin{equation}
\fair(x,z)=f(x+z)\ .
\end{equation}
Then
\begin{equation}
\lair(\fair)=z\fair \ , 
\ \lair(z,\d_z)(\fair)=x\fair\ .\label{airy is bispectral}\end{equation}

To generalize Wilson's formula to this case, introduce
a {\it Baker functional}, $\psiair$, defined on $ker \lair$, of 
the form 
\begin{equation} 
\Psi_{\text{Ai}}(f)=\sum_{i=0}^{r-1}k_i(x,z)\d^i(\fair)
\ .\label{airy functional}
\end{equation}
To specify its properties, it is useful to introduce the dual
description of
$\grad$ (cf. \cite{W1,KR}).  
As Wilson shows,  every point of $W\in\grad$ arises from a 
homogeneous, finite dimensional  space of finitely supported distributions 
in the
complex plane.  
That is, there should be complex numbers  $\lambda_1, \dots ,
\lambda_n$  and polynomials over  $\Bbb C$,  $\ell_1, \dots , \ell_n$, such
that 
\begin{equation}
W=\frac 1{q(z)}\{\ p(z)\in\Bbb C[z]\ |\ c_i(p)=0\ \}\ ,
\end{equation}
where
\begin{equation}
q(z)=\Pi_i(z-\lambda_i)
\end{equation}
and
\begin{equation}
c_i=\delta_{\lambda_i}\circ\ell(\d_z)\ .
\end{equation}

Now  define $\psiair$
by the following properties.
Let $C$ denote the span of  $c_1,...,c_n$. 

\vskip 12pt
\noindent{ Property 1a:}\ \ \  The  functions  $q(z) k_i(x,z)$ are  polynomial in
$z$.
\vskip 12pt
\noindent{ Property 2a:} 
For all  $f\in Ker(L)$,
$q_C(z)\psiair(f)$  is annihilated by all $c\in C$.
\vskip 12pt
\noindent{ Property 3a:} 
$\lim\limits_{z\to\infty}(k_0,...,k_{r-1})=(1,0,...,0)$.
\vskip 12pt
Fixing the operator $\lair$, set
\begin{equation}\label{airy vector}
\vec{k}_W=\left\lgroup\begin{matrix} k_0\\ \vdots
\\k_{r-1}\end{matrix}\right\rgroup.  
\end{equation}  
It is important to remark that the map $C\mapsto W$ is not
one-to-one.   The equivalence relation on $C$ induced by this map
is the one generated by
\begin{equation}\label{equivalence relation} 
C + \Bbb C\delta_{\lambda} \sim C \circ
(z -
\lambda)\ .
\end{equation}
In particular, the properties of $\psiair$ depend only on $W$ and 
not on the representative $C$.

Choose a space of conditions
$C$ defining $W$, and  set $k_i(x,z)=\frac 1{q_C(z)}\sum_jk_{i,j}(x)z^j$.
Consider the differential operator
\begin{equation}
K_{\text{Ai},C}=\sum_{i,j} k_{i,j}(x)\d^i\lair^j\ .
\end{equation}
Defining
\begin{equation}\label{airy vector f}
\vecair(x,z)=(\fair,\d(\fair),...,\d^{r-1}(\fair))\ ,
\end{equation}
one has 
\begin{equation}
K_{\text{Ai},C}(\fair)=q_C(z) \vecair\cdot\vec{k}_W\ .
\end{equation}
The asympotics of $\vec{k}_W$
imply that $K_{\text{Ai},C}$ is a monic
operator of order
$rn$.  Property 2a
implies that $K_{\text{Ai},C}$ annihilates
$c(\fair)$ for all the distributions $c\in C$ and all $f\in Ker(\lair)$.  
Thus
$K_{\text{Ai},C}$ is unique operator with the two properties just stated.   It then
follows (cf. \cite{KR}) that for any polynomial $p\in R_W$, the
pseudodifferential operator
$M_p=K_{\text{Ai},C}p(\lair)K_{\text{Ai},C}^{-1}$ is a differential operator, and 
\begin{equation} M_p(K_{\text{Ai},C}(\fair))=p(z)K_{\text{Ai},C}(\fair)
\end{equation} 
for all $f\in Ker(\lair)$.   Note that  $M_p$ could also have been obtained by
conjugating $p(\lair)$  by the 
monic $0$th order pseudodifferential operator $K_{\text{Ai},C}q(\lair)^{-1}$.  The latter
operator is independent of the space $C$ representing $W$, and is the analogue of 
the Sato operator in this theory.

Thus one has a rank $r$ commutative algebra of differential operators
\begin{equation}\label{airy algebra}
\cal R_{\text{Ai},W}=\{\ M_p\ |\ p\in\ R_W\ \}\ ,
\end{equation}
with an $r$-dimensional space of eigenfunctions $\psiair(f)$,
$f\in Ker(L)$.
To accomodate the usual normalization of spectral algebras, namely
that the subprincipal symbol should vanish, take $\rho$ of the
form $\rho(t)=t^r+O(t^{r-2})$.

The task is to obtain a formula for $\vec{k}_W$ in terms of the
matrices
$P$ and $Q$ corresponding to the point $W$.    
The following lemma expresses  property 2a in terms of 
covariant differentiation of  $\vec {k}_W$.  Set
\begin{equation}\label{first order airy system}   
\bair(x,z)=
\left[\begin{matrix} 
0&&\dots&&0&x+z\\   
1&\ddots&&&&a_1\\   
0&\ddots&&&\vdots&a_2\\
\vdots&\ddots&&&&\vdots\\   
&&&&0&a_{r-2}\\  
0&&\dots&0&1&0
\end{matrix}\right]\ ,
\end{equation}
where the $a_i$'s are the coefficients of $\rho$, and set
\begin{equation}\label{nabla}
\nabair=\frac{\d}{\d z}+\bair(x,z)\ .
\end{equation}
\begin{Lem}\label{covariant airy}  Let $c$ be a  distribution of the form  
$c=\delta_{\lambda}\circ p(\d_z)$. Let
$\vec{g}=
\left[\begin{matrix} g_0(x,z)\\
\vdots\\ g_{r-1}(x,z)\end{matrix}\right]$ be a vector of polynomials in $z$
with coefficients in $\C(x)$. Then
$c(\vecair\cdot\vec{g})=0$ for  all $f\in Ker(L)$ if
and only if $\vec{g}$ is annihilated by
$\delta_{\lambda}\circ p(\nabair)$.\end{Lem}
\begin{pf} Given $f\in Ker(L)$,
\begin{equation}
\d_z(\vecair)=\vecair\cdot\bair(x,z)\ .
\end{equation}  
Thus, for all $j$,
\begin{equation}
\delta_{\lambda}\circ\d_z^j(\vecair\cdot\vec{g})=
\vecair(x,\lambda)\cdot(\d_z+\bair(x,\lambda))^j(\vec{g})\ .
\end{equation} This proves the lemma, since there is no differential equation
of order less than $r$ satisfied by $f(x+\lambda)$ for all $f\in Ker(L)$.
\end{pf}

Now consider the involutions on each $C_n$ defined by
\begin{align}\label{airy involution}
  C_n&\overset{\beta_{\text{Ai}}}\too C_n\notag\\
((P,Q))&\mapsto((\hat{P},\hat{Q}))\ ,\ \ \ \text{where}\notag\\
\hat{P}=P^T;\ \ \ &\hat{Q}=\rho(P^T)-Q^T\ .
\end{align}
\begin{Thm}\label{airy formula}
Let  $W \in \grad$  correspond to a point  $((P, Q)) \in C_n$.  Let
\begin{equation} [P, Q] = I - w_1w^T_2
\end{equation} where  $w_1$  and  $w_2$  are column vectors.  
For  $j = 0, \dots ,
r - 1$, let
\begin{equation}
\rho_j(t) = t^{r-1-j} - \sum^{r-2}_{i=j+1} a_it^{i-1-j} .
\end{equation} Then the components of $\vec{k}_W$ are
\begin{equation}\label{k formula}
 k_j(x, z) = \delta_{0,j} - w^T_2(zI -
Q)^{-1}\rho_{j}(P)(x -\hat{Q}^T)^{-1}w_1 \ .
\end{equation}
\end{Thm}
\begin{pf}It suffices to consider the generic case, in which  $W$  is defined by
conditions  $c_i = \delta_{\lambda_i}
\circ(\partial - \alpha_i)$, $i = 1, \dots , n$, for distinct $\lambda$'s.
Set
\begin{equation}\label{gamma}
\gamma_i = \alpha_i - \sum_{j \ne i} \frac{1}{\lambda_i -
\lambda_j} .
\end{equation} Then  $W$ corresponds under
Wilson's theorem to the Calogero-Moser
pair
\begin{equation} Q =
\left[
\begin{matrix}
\lambda_1 &&0\\ &\ddots&  \\ 0 & & \lambda_n
\end{matrix}
\right]
\end{equation}
\begin{equation} P = 
\left[
\begin{matrix}
\gamma_1 & \frac{1}{\lambda_1 - \lambda_2} & \dots &
\frac{1}{\lambda_1 - \lambda_n} \\
\frac{1}{\lambda_2 - \lambda_1} & & & \frac{1}{\lambda_2 -
\lambda_n} \\
\vdots&&\ddots&\vdots \\
\frac{1}{\lambda_n - \lambda_1} & \dots & & \gamma_n
\end{matrix}
\right]\ .
\end{equation}

Let  $(e_j)_{j=1,\dots , r}$  be the standard basis for  $\Bbb C^r$.  Then
\begin{equation}\label{k form}
\vec{k}_W= e_1 + \sum^n_{i=1} \frac{\vec{v}_i(x)}{(z - \lambda_i)}\ ,
\end{equation} for some vectors  $\vec{v}_i(x) \in \Bbb C^r$.  Applying 
$\delta_{\lambda_i} \circ \left(\d_z + \bair(x, 
\lambda_i) - \alpha_i\right)$  to  $q(z)\vec{k}_W$, one obtains the set of
equations  
\begin{align} 0 &= \bair(x , \lambda_i)\vec{v}_i(x) \prod_{\ell\ne i}
(\lambda_i -
\lambda_{\ell}) + \notag \\ &\prod_{\ell\ne i} (\lambda_i - \lambda_{\ell})e_1 +
\sum^n_{m=1}
\vec{v}_m(x) \sum_{\ell\ne m} \prod_{s\ne \ell, m} (\lambda_i -
\lambda_s) \notag \\ &- \alpha_i \vec{v}_i(x) \prod_{\ell\ne i} (\lambda_i -
\lambda_{\ell}) ,
\end{align}
$i = 1, \dots , n$.  Dividing by $\prod_{\ell\ne i} (\lambda_i -
\lambda_{\ell})$ and using \eqref{gamma},
\begin{equation}\label{separate equations} e_1 = -\bair(x ,
\lambda_i)\vec{v}_i(x) +
\gamma_i\vec{v}_i(x) -
\sum_{\ell\ne i} \frac{\vec{v}_{\ell}(x)}{\lambda_i -
\lambda_{\ell}} .
\end{equation} 
Let 
$(u_i)_{i=1,\dots , n}$  be the standard basis for  $\Bbb C^n$, and let
\begin{equation} v(x) = \sum u_i \otimes \vec{v}_i(x) \in \Bbb C^n \otimes \Bbb C^r .
\end{equation} Let  $w = \sum u_i$.   Write $\bair(x , \lambda_i) =
\bair(x,0) +
\lambda_i\Delta$, where
\begin{equation}
\Delta = 
\left(
\begin{matrix} 0 & \cdots & 0 & 1 \\ 0&&0& 0 \\
\vdots&\cdots&\vdots& \vdots \\ 0 &\cdots&0& 0
\end{matrix}
\right) .
\end{equation} Then  \eqref{separate equations} is encoded as a single equation
\begin{equation} w \otimes e_1 = (-I \otimes B(x) + P \otimes I - Q \otimes
\Delta)v(x) .
\end{equation} Altogether, \eqref{k form} becomes
\begin{equation}
\vec{k}_W = e_1 - (w^T \otimes I) \circ ((zI - Q)^{-1} \otimes I) \circ A^{-1}
(w
\otimes e_1) ,
\end{equation} where
\begin{equation} A = I \otimes \bair(x,0) - P \otimes I + Q \otimes \Delta \in
End(\Bbb C^n \otimes \Bbb C^r) .
\end{equation} Thinking of  $\Bbb C^n \otimes \Bbb C^r$  as  $rn$-tuples in blocks
of length  $n$,
\begin{equation} w \otimes e_1 = 
\left[
\begin{matrix} 1 \\
\vdots \\ 1 \\ 0 \\
\vdots \\ 0
\end{matrix}
\right] 
\begin{matrix}
\biggr\}n \\
\phantom{} \\
\phantom{} \\
\phantom{}
\end{matrix}
\end{equation} and
\begin{equation}\label{airy matrix} A = 
\left[
\begin{matrix} -P & 0 & \cdots & 0 & xI + Q \\ I & -P &&& a_1I \\ 0 & I & \ddots &&
\vdots \\
\vdots & 0 &&& a_{r-1}I \\ 0 & 0 & 0 & I & -P
\end{matrix}
\right] .
\end{equation} One checks quite easily that
\begin{equation} A
\left[
\begin{matrix}
\rho_0(P) \\
\vdots \\
\vdots \\
\rho_{r-1}(P)
\end{matrix}
\right] =
\left[
\begin{matrix} xI + Q - \rho(P) \\ 0 \\
\vdots \\ 0
\end{matrix}
\right] .
\end{equation} Then
\begin{align} &(w^T \otimes I) \circ ((zI - Q)^{-1} \otimes I) \circ A^{-1}(w
\otimes e_1) \notag \\ &= \left[
\begin{matrix} w^T&&0 \\ &\ddots& \\ 0 && w^T
\end{matrix}
\right]
\left[
\begin{matrix} (zI - Q)^{-1}&&0 \\ &\ddots&  \\ 0 && (zI - Q)^{-1}
\end{matrix}
\right]
\left[
\begin{matrix}
\rho_0(P) \\
\vdots \\
\vdots \\
\rho_{r-1}(P)
\end{matrix}
\right] (xI + Q - \rho(P))^{-1}w . \notag \\ &= \left[
\begin{matrix} w^T(zI - Q)^{-1}\rho_0(P)(xI-\hat{Q}^T )^{-1}w \\
\vdots \\ w^T(zI - Q)^{-1}\rho_{r-1}(P)(xI-\hat{Q}^T )^{-1}w
\end{matrix}
\right] .
\end{align} This proves the theorem, since
\begin{equation} [P,Q] = I - ww^T .
\end{equation}
\end{pf}

\begin{Cor}  For any generalized Airy operator $\lair$, and for all
$W\in\grad$, 
\begin{equation}\label{xz symmetry}
\vec{k}_W(z,x)=\vec{k}_{\beta_{\text{Ai}}(W)}(x,z)\ .
\end{equation}
In particular, the algebra  $\cal R_{\text{Ai},W}$, \eqref{airy algebra}, is
bispectral, with an
$r$-dimensional space of joint eigenfunctions $\psiair(f)$, $f\in
Ker(\lair)$.
\end{Cor}
\begin{pf} Formula \eqref{xz symmetry} follows immediately from \eqref{k formula},
for if $[P,Q]=I-w_1w_2^T$,  then 
\begin{equation}
[\hat{P},\hat{Q}]=[P^T,\rho(P^T)-Q^T]=[P,Q]^T=I-w_2w_1^T\ .
\end{equation}
The rest of the corollary is immediate.
\end{pf}

\section{The Bessel Case}

The Bessel case
works in much the same way.   Consider again a polynomial $\rho(t)$, now 
normalized so that $a_{r-1}=\binom r2$.
Set
\begin{equation}
\lbess=x^{-r}\rho(D)\ \ \ \ \ \ \ \text{(Bessel)}\
,\label{besselop}\end{equation} where
\begin{equation}
\d=\frac d{dx}\ ,\ D=x\d\ .
\end{equation}
Consider $Ker(\lbess-1)$, which should now be thought of as 
a sheaf rather than a space.    For $f\in Ker(\lbess-1)$, define
\begin{equation}
\fbess(x,z)=f(xz)\ .
\end{equation}
Then
\begin{equation}
\lbess(\fbess)=z^r\fbess  \ ,
\ \lbess(z,\d_z)(\fbess)=x^r\fbess\ .\label{bessel is
bispectral}\end{equation}
Assume now that the matrix $Q$ is invertible.   
Define, as
the analogue of \eqref{airy vector f},
\begin{equation}\label{bessel vector f}
\vecbess(x,z)=(\fbess,D(\fbess),...,D^{r-1}(\fbess))\ .
\end{equation}
Then
\begin{equation}\label{bessel exchange}
D_z(\vecbess)=\vecbess\bbes(x^r,z^r)\ ,
\end{equation} 
where
\begin{equation}\label{first order bessel system}    
\bbes(x,u)=
\left[\begin{matrix}  0&&\dots&&0&a_0+xu\\    1&\ddots&&&&a_1\\   
0&\ddots&&&\vdots&a_2\\
\vdots&\ddots&&&&\vdots\\    &&&&0&a_{r-2}\\   
0&&\dots&0&1&a_{r-1}
\end{matrix}\right]\ .
\end{equation}
Accordingly, one expects a  Baker functional of the form
\begin{equation}\label{bessel baker functional}
\psibess(f)=\vecbess\cdot\vec{k}_W(x^r,z^r)\ \ ,\ f\in Ker(\lbess-1)\ .
\end{equation}

To state the properties of $\psibess$,   introduce the functions
\begin{equation}
\mu(x,z)=(x^r,z^r)\ ;\ \nu(z)=z^r\ .
\end{equation}
Denote by $\nu^*$ the action of $\nu$ on the space of finitely supported
distributions in $\C^*$.   Let
\begin{equation}
\nabbess=D_z+ \frac 1r \bbes(x,z)\ .
\end{equation}
Given a distribution $c=\delta_{\lambda}\circ p(D)$, define
\begin{equation}
c_{\nabbess}=\delta_{\lambda}\circ p(\nabbess)\ ,
\end{equation}
acting on vector valued functions of $z$.

\begin{Lem}\label{covariant bessel}  Let $c$ be a  distribution of the form  
$c=\delta_{\lambda}\circ p(D)$. Let
$\vec{g}=
\left[\begin{matrix} g_0(x,z)\\
\vdots\\ g_{r-1}(x,z)\end{matrix}\right]$ be a vector of polynomials in $z$ with
coefficients in $\C(x)$. Then
$c(\vecbess\cdot\mu^*(\vec{g}))=0$ for  all $f\in Ker(\lbess-1)$ if and only if
$\nu^*(c)_{\nabbess}(\vec{g})=0$.\end{Lem}
\begin{pf} By virtue of the identity
\begin{equation}
D\circ \nu^*=\nu^*\circ r D\ ,
\end{equation}
one has
\begin{equation}
\nu^*(c)=\delta_{\lambda^r}\circ p(rD)\ .
\end{equation}
By \eqref{bessel exchange},
\begin{align}
c\circ\vecbess\cdot\mu^*(\vec{g})&
 =\delta_{\lambda}\circ\vecbess\cdot p(D+\bbes(x^r,z^r))\circ
\mu^*(\vec{g})\notag\\
&=\vecbess(x,\lambda)\cdot
\delta_{\lambda}\circ \mu^*\circ p(rD+\bbes(x,z))
(\vec{g})\notag\\
&=\vecbess(x,\lambda)\cdot\nu_x^*\circ
\delta_{\lambda^r}\circ p(r\nabbess)
(\vec{g})\notag\\
&=\vecbess(x,\lambda)\cdot\nu_x^*\circ \nu^*(c)_{\nabbess} (\vec{g})\ ,
\end{align}
where $\nu_x$ is $\nu$ acting in the $x$-variable.  The lemma now follows 
as in lemma \ref{covariant airy}.\end{pf}

In light of the preceding lemma, it makes sense to impose the following
properties on $\psibess$.
\vskip 12pt

\noindent{ Property 1b:}\ \ \  The  functions  $q_C(z) k_i(x,z)$ are 
polynomial in
$z$.
\vskip 12pt
\noindent{ Property 2b:}   Let $C'$ be any space of distributions such that
$\nu^*(C')=C$.   Then  
for all $f\in Ker(\lbess-1)$, $q_C(z)\psibess(f)$ is annihilated by all 
$c\in C'$.
\vskip 12pt
\noindent{ Property 3b:} 
$\lim\limits_{z\to\infty}\vec{k}_W=e_1$.
\vskip 12pt
As in the Airy case, one reconstructs a differential
operator $K_{\text{Be},C}$, but now
\begin{equation} K_{\text{Be},C}(\fbess)=q_C(z^r)
\vecbess\cdot\vec{k}_W(x^r,z^r)\ .
\end{equation}
Then for any polynomial $p\in R_W$, the
pseudodifferential operator
$M_p=K_{\text{Be},C}p(\lbess)K_{\text{Be},C}^{-1}$ is a differential operator, and 
\begin{equation} M_p(K_{\text{Be},C}(\fbess))=p(z^r)K_{\text{Be},C}(\fbess)
\end{equation}  for all $f\in Ker(\lbess-1)$.   Define 
$\cal R_{\text{Be},W}$ to be the algebra of the all the $M_p$'s.

Everything now proceeds as before. 
Assume  that the matrix $Q$ is invertible.  We have $n$ distributions
$c_i = \delta_{\lambda_i}
\circ(\partial_z - \alpha_i)$.     Note that 
$\delta_{\lambda_i}
\circ(\partial_z - \alpha_i)=\frac 1{\lambda_i} \delta_{\lambda_i}
\circ(D_z - \lambda_i\alpha_i)$.   Thus, according to 
lemma \ref{covariant bessel},  property 2b imposes the $n$ conditions
\begin{equation}
0=\delta_{\lambda_i}
\circ(\partial_z +\frac 1{r \lambda_i}\bbes(x,\lambda_i)-
\alpha_i)(q_C(z)\vec{k}_W(x,z))\ .
\end{equation}
Setting
\begin{equation}\label{bessel k form}
\vec{k}_W = e_1 + \sum^n_{i=1} \frac{\vec{v}_i(x)}{(z - \lambda_i)} \ ,
\end{equation} 
one now finds 
\begin{equation}\label{separate bessel
equations} e_1 = -\frac 1{r\lambda_i} \bbes(x,\lambda_i)\vec{v}_i(x) +
\gamma_i\vec{v}_i(x) -
\sum_{\ell\ne i} \frac{\vec{v}_{\ell}(x)}{\lambda_i -
\lambda_{\ell}} .
\end{equation} 
This time,
\begin{equation}
\frac 1{ru}\bbes(x,u) = \frac 1{ru}\Delta_1 + \frac xr\Delta_2\ ,
\end{equation}
where
\begin{align}
\Delta_1&= 
\left[
\begin{matrix} 0 & \cdots && 0 & a_0 \\ 
1&\ddots&&\vdots& \vdots \\
0&\ddots&&0&\vdots\\
&\ddots&&&\\
0 &\cdots&0&1& a_{r-1}
\end{matrix}
\right] \ ,\\ 
&\notag\\
\Delta_2&=\left[
\begin{matrix} 
0&\cdots&0&1\\
\vdots&&\vdots&0\\
\vdots&&\vdots&\vdots\\
0&\cdots&0&0
\end{matrix}
\right] \ .
\end{align}
Thus  
\begin{equation} w \otimes e_1 = 
-(\frac xr I \otimes \Delta_2 - P \otimes I + \frac 1r  Q^{-1} \otimes
\Delta_1)v(x) \ .
\end{equation} 
Then
\begin{equation}
\vec{k}_W = e_1 - (w^T \otimes I) \circ ((zI - Q)^{-1} \otimes I) \circ A^{-1}
(w
\otimes e_1) ,
\end{equation} where  $A$ is now given in block matrix form by
\begin{equation}\label{bessel matrix} 
A = 
\left[
\begin{matrix} 
-P & 0 & \cdots & 0 & \frac xr I + \frac{a_0}r Q^{-1} \\ 
\frac 1r Q^{-1} & \ddots &&\vdots& \frac{a_1}r Q^{-1} \\ 
0 &\ddots & &0&\vdots \\
\vdots & \ddots &&-P& \frac{a_{r-2}}r Q^{-1} \\ 
0 & \cdots & 0 & \frac 1r Q^{-1} & -P+\frac{a_{r-1}}r Q^{-1}
\end{matrix}
\right] .
\end{equation} 
One obtains the following result.

\begin{Thm}\label{bessel formula} 
Let  $W \in \grad$  correspond to a point  $((P, Q))
\in C_n$.  Let
\begin{equation} [P, Q] = I - w_1w^T_2
\end{equation} where  $w_1$  and  $w_2$  are column vectors. 
Writing the  $r$th order Bessel operator $\lbess$ in the form 
$\lbess=x^{-r}\rho(D)$, let  
\begin{equation}
\rho_j(t) = t^{r-1-j} - \sum^{r-1}_{i=j+1} a_it^{i-1-j} 
\end{equation} for  $j = 0, \dots , r - 1$. 
Then the components of $\vec{k}_W$
are
\begin{equation}\label{bessel k formula}
 k_{j}(x, z) = \delta_{0,j} - 
r w^T_2(zI - Q)^{-1}\rho_{j}(r QP)(x
-\hat{Q}^T)^{-1}w_1 ,
\end{equation}
where
\begin{equation}
\hat{Q}=(Q^{-1}\rho(r QP))^T\ .
\end{equation}
\end{Thm}
\begin{pf}
One checks now that with $A$ given by \eqref{bessel matrix},
\begin{equation} A
\left[
\begin{matrix}
r \rho_0(r QP) \\
\vdots \\
r \rho_{r-1}(r QP)
\end{matrix}
\right] =
\left[
\begin{matrix} xI - Q^{-1} \rho(r QP) \\ 0 \\
\vdots \\ 0
\end{matrix}
\right] .
\end{equation}
The result then follows as in theorem \ref{airy formula}.\end{pf}

Theorem \ref{bessel formula} suggests the definition
\begin{align}
\hat{P}&=\hat{Q}^{-1}P^TQ^T\notag\\
&=(QP\rho(QP)^{-1} Q)^T \  .
\end{align}
Then
\begin{align}
\hat{\hat{Q}}&=(\hat{Q}^{-1}\rho(r \hat{Q}\hat{P}))^T\notag\\
&=(\hat{Q}^{-1}\rho(r P^TQ^T))^T\notag\\
&=\rho(r QP)\rho(r QP)^{-1} Q=Q\ ,
\end{align}
and
\begin{align}
\hat{\hat{P}}&=\hat{\hat{Q}}^{-1}\hat{P}^T \hat{Q}^T\notag\\
&=Q^{-1}QP(\hat{Q}^{-1})^T\hat{Q}^T=P\ .
\end{align}
Moreover,
\begin{align}
[\hat{P},\hat{Q}]&=[\hat{Q}^{-1}P^TQ^T,\hat{Q}]\notag\\
&=\hat{Q}^{-1}P^TQ^T\hat{Q}-P^TQ^T\notag\\
&=Q^T(\rho(rQP)^{-1})^T(QP)^T\rho(rQP)^T(Q^{-1})^T-P^TQ^T\notag\\
&=Q^TP^T-P^TQ^T=[P,Q]^T\ .
\end{align}
So again one has an involution (densely defined) on $C_n$,
\begin{align}\label{bessel involution}
  C_n&\overset{\beta_{\text{Be}}}{--\rightarrow} C_n\notag\\
(P,Q)&\mapsto(\hat{P},\hat{Q})\ ,\ \ \ \text{where}\notag\\
\hat{P}=\hat{Q}^{-1}P^TQ^T\ ;
\ \ \ &\hat{Q}=
(Q^{-1}\rho(r QP))^T
\ .
\end{align}
\begin{Cor}  For any generalized Bessel operator $\lbess$,
\begin{equation}\label{bessel xz symmetry}
\vec{k}_W(z,x)=\vec{k}_{\beta_{\text{Be}}(W)}(x,z)\ .
\end{equation} In particular, the algebra  $\cal R_{\text{Be},W}$ is bispectral,
with a rank $r$ joint eigensheaf
$\psibess(f)$,
$f\in Ker(\lbess-1)$.
\end{Cor}
\section{Dynamics}

It is well-known that the Calogero-Moser particle system is a completely
integrable hamiltonian system on the symplectic manifold $C_n$.  The
symplectic form is
\begin{equation}
\omega=\tr(dP dQ)\ ,
\end{equation}
and the hamiltonians are $h_n=\tr(P^n)$ (cf. \cite{KKS}).

It is pleasing then that the Airy and Bessel involutions are antisymplectic on
each
$C_n$.    In the Airy case,
\begin{equation}
\hat{P}=P^T;\ \ \ \hat{Q}=\rho(P^T)-Q^T\ ,
\end{equation} this follows from the fact that
\begin{equation}
\tr(dP d P^n)=\sum_{i+j=n-1}\tr(dP P^i dP P^j) =0 \ .
\end{equation} (The basic trace identity for form-valued matrices is
\begin{equation}
\tr(XY)=(-1)^{\text{deg}(X)\text{deg}(Y)}\tr(YX)\ .)
\end{equation}

The Bessel case is slightly more involved.
\begin{Prop}   Let $\sigma$ be a polynomial over $\C$.   Let
\begin{equation}
\hat{Q}= (Q^{-1}\sigma(QP))^T\ ;\ \hat{P}=\hat{Q}^{-1}P^TQ^T\ .
\end{equation} Then
\begin{equation}
\tr(d\hp d\hq)=-\tr(dP dQ)\ .
\end{equation}
\end{Prop}
\begin{pf} Let $\tq=Q^{-1}\sigma(QP)$, $\tp=QP\sigma(QP)^{-1} Q=\hp^T$.  Then
\begin{equation}
\tr(d\hp d\hq)=\tr(d \tp d\tq)\ .
\end{equation} Set $R=QP$.  Then
\begin{equation}
\tr(dP dQ)=\dtr (Q^{-1} R dQ)\ ,
\end{equation} while
\begin{align}
\tr(d\tp d\tq) &= \dtr(R\sigma^{-1} Qd(Q^{-1}\sigma)) \nonumber \\ &=
\dtr(R\sigma^{-1}d\sigma - R\sigma^{-1}dQ Q^{-1}\sigma)
\nonumber \\ &= \dtr(R\sigma^{-1}d\sigma - Q^{-1} R dQ) \ .
\end{align} So it must be proved that
\begin{equation}
\dtr(R\sigma^{-1}d\sigma) = 0 \ .\label{to prove}
\end{equation} If  $\sigma = \sigma_1\sigma_2$  and  $\sigma_2$  
commutes with 
$R$, then
\begin{align}
\dtr(R\sigma^{-1}d\sigma) &= 
\dtr(R\sigma^{-1}_2\sigma^{-1}_1(d\sigma_1\sigma_2 + \sigma_1 d\sigma_2))
\nonumber \\ &= \dtr(R\sigma^{-1}_1d\sigma_1) + 
\dtr(R\sigma^{-1}_2d\sigma_2)\ .  \label{additivity}
\end{align} Also,
\begin{equation}
\dtr(\sigma^{-1}d\sigma) = -\tr(\sigma^{-1}d\sigma
\sigma^{-1}d\sigma) = 0 .
\end{equation} This last identity implies that one can replace  $R$  by 
$R-(const)  I$  in \eqref{to prove}, and reduce to the case that  $\sigma =
R\sigma_1(R)$, $\sigma_1$  polynomial.

Then by \eqref{additivity}
\begin{align}
\dtr(R\sigma^{-1}d\sigma) &= \dtr(dR + R\sigma^{-1}_1 d\sigma_1) \nonumber \\
&= \dtr(R\sigma^{-1}_1 d\sigma_1) \ . 
\end{align} Since the result is obvious when  $\sigma$  is a constant, the
proposition follows by induction on the degree of  $\sigma$.
\end{pf}

Now  introduce time dependence into the Baker functionals of the previous
sections in a manner generalizing the standard procedure in the rank-one case
\cite{SW}.
Fix a positive integer $m$.   The time-dependent Baker function in the rank one
case is the function
\begin{equation} e^{xz+tz^m}p(x,z,t)
\end{equation} satisfying properties 1,2 and 3 of section \ref{airy case} with a
fixed space of conditions $C$.   On the other hand, one may introduce
time-dependence into the conditions by defining
$C_t=C\circ e^{tz^m}$.   Then the  function
$e^{xz}p(x,z,t)$ satisfies properties 1,2 and 3 for the variable conditions
$C_t$.  

If $c=\delta_{\lambda}\circ(\d-\alpha)$, then
$c\circ e^{t z^m}=\delta_{\lambda}\circ(\d+t m \lambda^{m-1}-\alpha)$. In
other words, the flow $C_t$ is seen on the level Calogero-Moser pairs as 
\begin{equation}\label{cm flows} Q_t=Q_0\ ;\ P_t=P_0-tmQ^{m-1}\ .
\end{equation} This is the flow  of the completely integrable hamiltonion
\begin{equation} h_m=\tr(Q^m)\ .
\end{equation} This hamiltonian is the Calogero-Moser hamiltonian with the
roles of $P$ and $Q$ reversed.  Finally, the Baker function gives rise in a
standard way to a solution of the KP-hierarchy, with poles in $x$ the same as
those of the Baker function.  Thus, Wilson's formula makes it immediately
clear that the poles in
$x$ of such a KP-solution move as a Calogero-Moser particle system (cf.
\cite{Kr,Sh,Ka}).

To carry this over to the Airy and Bessel cases, define $\vec{k}_{W,t}$ to
be the vector
\eqref{airy vector}, for the variable space of conditions $C_t$.   We are led
to solutions of a subhierarchy of the KP-heirarchy, in the following way. Let
$\kair$ be the monic $0^{th}$-order pseudodifferential operator such that
\begin{equation}
\kair\d^r\kair^{-1}=\lair\ .
\end{equation} Let $\tilde{K}_t$ be the monic $0^{th}$-order
pseudodifferential operator such that
\begin{equation}
\vecair\cdot\vec{k}_{W,t}=\tilde{K}_t(\fair)\ .
\end{equation} Then the argument in \cite{SW} shows that
\begin{equation}\label{pre dressing}
\d_t(\tilde{K})\tilde{K}^{-1}+(\tilde{K}\lair^m\tilde{K}^{-1})_-=0\ .
\end{equation} Now let
\begin{equation} K_t=\tilde{K}_t\kair\ .
\end{equation}
Then
\begin{equation}\label{dressing}
\d_t(K)K^{-1}+(K\d^{rm}K^{-1})_-=0\ .
\end{equation} It now follows as in \cite{SW} that the operator
\begin{equation} M_t=K_t\d K_t^{-1}
\end{equation} satisfies the $rm^{th}$ term of the KP-hierarchy,
\begin{equation}
\d_t(M)=[M^{rm}_+,M]\ .
\end{equation} From formula \eqref{k formula}, $M$ has coefficients in
$\C(t)[x,\frac 1{\tau_t(x)}]$, where
\begin{align}
\tau_t(x)&=Det(x-\hat{Q}_t)\ ,\\
\hat{Q}_t&=(\rho(P_t)-Q_t)^T\notag\\ &=(\rho(P_0-tmQ_0^{m-1})-Q_0)^T\ .
\end{align} Thus we have constructed a solution $M_t$, of the $rm^{th}$ term
of the KP-hierarchy, whose poles in $x$ move according to the completely
integrable hamiltonian $\tr(\hq^m)$ on the Calogero-Moser phase space.

In the Bessel case, exactly the same analysis holds, except that 
$M$ has coefficients in
$\C(t)[x,\frac 1x,\frac 1{\tau_t(x^r)}]$. Thus $M_t$ has a fixed pole at
$x=0$, with the motion of the remaining poles being governed by the
hamiltonian 
$\tr(\hq^m)$.

With several particles it becomes quite cumbersome to write out these
hamiltonians explicitly.    The lowest rank cases are
\begin{align}
\lair=\d^2-x\ &;\ \rho(t)=t^2\ ,\\
\lbess=\d^2-x^{-2}\ &;\ \rho(t)=t^2-t-1\ .
\end{align} 
Note that  the
first  hamiltonian, $h_1=\tr(\hq)$, is already non-linear.  In the rank-two
Bessel case with one particle, for instance, 
$h_1(\lambda,\gamma)=\lambda^{-1}\rho(\lambda\gamma)$.  This gives the
equations of motion
\begin{align}
\dot{\lambda}&=8\gamma\lambda-2\\
\dot{\gamma}&=-4\gamma^2-\lambda^{-2}\ .
\end{align} These equations are solved by applying the Bessel involution and
changing $t$ to
$-t$,  i.e. by setting
\begin{equation}
\hat{\lambda}=c_1\ ;\ \hat{\gamma}=c_2+t\ .
\end{equation} After some calculation, 
\begin{equation}
\lambda(t)= 4 c_1 t^2  +  ( 8 c_1 c_2-2) t -\frac 1{c_1} - 2 c_2  + 4 c_1
c_2^2\ .
\end{equation} With two particles and a second order Airy operator, the first
hamiltonian is
\begin{equation} {\gamma_1^2} + {\gamma_2^2} -\lambda_1 - \lambda_2 -  {2\over
{{{\left( -\lambda_1 + \lambda_2 \right) }^2}}}\ .
\end{equation}
With two particles and a second order Bessel operator, the first hamiltonian
is
\begin{equation}
 -{{\lambda_1 + \lambda_2}\over {\lambda_1\,\lambda_2}} - \gamma_1 + 
\lambda_1\,{\gamma_1^2} - 
\gamma_2 +
\lambda_2\,{\gamma_2^2} + 
  {{2\,\left( \lambda_1\,\gamma_1 - 
\lambda_2\,\gamma_2 \right) }\over {-\lambda_1 + \lambda_2}}\ .
\end{equation}

\end{document}